\documentclass[a4paper,fleqn]{article}
\usepackage[ansinew]{inputenc}
\usepackage{amsmath,amssymb}
\usepackage{theorem}
\usepackage{color}
\newtheorem{theorem}{Theorem}
\newtheorem{lem}[theorem]{Lemma}

\newenvironment{beweis}{\textbf{Proof:}}{\hfill$\Box$\\}

\newtheorem{defi}[theorem]{Definition}
\newtheorem{kor}[theorem]{Corollary}
\newtheorem{beisp}[theorem]{Example}
\newtheorem{rem}[theorem]{Remark}
\newtheorem{prop}[theorem]{Proposition}
\newtheorem{alg}[theorem]{Algorithm}
\newtheorem{opt}[theorem]{Optimization Problem}

\theorembodyfont{\upshape}
\numberwithin{theorem}{section}
\numberwithin{equation}{section}

\newcommand{\R}{\mathbb{R}}

\newcommand{\N}{\mathbb{N}}

\title{The Weber Problem in Logistic and Services Networks
under Congestion}

\author{Vanessa Lange
\\
{Technical University of Mittelhessen}\\
 Department of  Mathematics,
		Natural Sciences, and Data Sciences,\\
Wilhelm-Leuschner-Strasse 13, 61169 Friedberg, Germany\\		
	vanessa.lange@mnd.thm.de	
	\and 
Hans Daduna
\thanks{Corresponding author}
\\
{Universit{\"a}t Hamburg, }
Department of Mathematics\\
Bundesstrasse 55,
20146 Hamburg,
Germany\\
 daduna@math.uni-hamburg.de
}

\begin{document}

\maketitle
\begin{abstract}
We investigate a location-allocation-routing 
problem where trucks deliver goods from a
   central production facility to a set of warehouses with fixed locations and known demands. Due to limited capacities congestion occurs and results in queueing problems.
The location of the center is determined to maximize the utilization of the given resources  (measured in throughput), and the minimal number of trucks is determined to satisfy the overall demand generated by the warehouses.
Main results for this integrated decision problem on strategic and tactical/operational level are:
(i) The location decision is reduced to a standard Weber problem with weighted distances.    
(ii) The joint decision for location and fleet size is separable.
(iii) The location of the center is robust against perturbations of  several  system parameters on the operational/tactical level.
Additionally, we consider minimization of  travel times as optimization target.
\end{abstract}
\underline{Keywords}: Weber problem, facility location, Gordon-Newell network, throughput optimization,
steady state analysis, travel times,  fleet size optimization,
robustness.

\section{Introduction}\label{sect:intro}
We investigate a generalized Weber problem in a dynamic and stochastic environment. The model combines a classical planar location problem with optimizing the size of the fleet of trucks in a logistic and services network
and describes interacting production  and distribution facilities with logistic components.\\
This  integrated 
problem is related to  location-allocation problems, 
location-routing problems (LRPs), location-inventory problems (LIPs), and  transporta\-tion-location-allocation problems
under random influences.
A key property of our model is to integrate strategic (facility location) and tactical/operational (allocation, scheduling) aspects of decision making. Such integrated logistics-location models occur e.g. in supply chain planning and operation, for  surveys see
\cite{melo;nickel;saldanha-da-gama:09} and \cite{heckmann;nickel:19}.
Difficulties arising with similar integration procedures in combined location-routing problems are described in  \cite{min;jayaraman;srivastava:98}[p.10].\\
Our model comprises a set of warehouses with known locations, a single production center which produces commodities which are demanded by the warehouses, and a set of trucks which transport the   commodities from the production center to the warehouses. 
The model captures additionally congestion which emerges at service stations with limited capacity for loading of trucks at the center and for unloading at the warehouses.
Optimality of the system is defined 
with respect to \textbf{maximization} of the warehouses' throughput (which generates revenue) and \textbf{minimization} of the number of trucks (which generate costs) \textbf{over time } and \textbf{in a stationary system}. {Our main contributions are}:\\
For given demands from the  warehouses and incorporating 
the consequences of congestion at the center and the warehouses
\\
\textbf{(i)} we find  an optimal location for the central production facility, and \\ 
\textbf{(ii)}  determine the needed transportation capacity,
and\\ 
\textbf{(iii)}  
solve jointly both optimization problems in a unified model, to combine strategic decision making  (for location of the center) with decisions for  tactical and operational issues (transportation capacity and scheduling rules for routing trucks), and\\
\textbf{(iv)} we demonstrate that the solution of the optimization problem reveals  important  robustness of that solution against changes or perturbations  of several parameters  on the tactical/operational level (insensitivity).\\
\underline{Summarizing}: We introduce a model which is tailored to integrate several hierarchical levels of decision making in complex systems. As indicated in the literature, neglecting these inter-dependencies often generates sub-optimal solutions, see the comments on literature in Section \ref{sect:Literature}, especially concerning 
integration of strategic and tactical/operational aspects of planning on p. \pageref{page:Integration}.

\textbf{Structure of the paper.} 
In Section \ref{sect:Literature} we review related literature. 
Section \ref{sect:Problem} presents an overview of the problem setting and connections to  
investigations of related problem settings.
Details are provided in Section \ref{sect:DetailsProblem}. 
Our main findings are described in Section \ref{sect:mainresults}. In Section \ref{sect:reliable stations} we analyze the system in full detail.
In Section \ref{sect:NumberOfTrucks} we complement our structural results by an algorithm to determine the minimal number of trucks needed to satisfy the overall demand at the warehouses.
In Section \ref{sect:DiscussionModel} we discuss the assumptions used in Section \ref{sect:reliable stations} and we indicate that in many cases  the results
hold in more general settings.
Section \ref{sect:MiniPassageTimes} is devoted to determine the optimal location of the center for minimizing the round-trip times of trucks.
Section \ref{sect:NumericManagerSensitive} contains numerical experiments to underpin the structural results obtained sofar and a discussion of robustnes properties of our results.
In Section \ref{sect:ManagerialDecision} we discuss the impact of our results on the consequences of intertwining, respectively separating  strategic (location) and 
tactical/operational (routing and scheduling) decisions.
We conclude by discussing generalizations and extensions of the location problem.
Necessary prerequisits from queueing network theory are provided in Appendix \ref{sect:fundamentals}.
Proofs are postponed to Appendix \ref{sect:proofs}.
Details of numerical experiments are presented in Section \ref{sect:AddNumeric}.\\
\textbf{Conventions:} 
$\N=\{0,1,2,\dots\}$, $\N_+=\{1,2,\dots\}$.
Empty sums are $0$, empty products are $1$. We set
$0/0=0$. For any set $A$, $\mathcal{P}(A)$ is the set of  subsets of $A$.
Increasing means non-decreasing and decreasing means non-increasing.

\section{Literature review}\label{sect:Literature}
\underline{Stochastic location models with congestion}, 
 i.e.
location problems in the context of service processes under stochastic influences have been investigated since the 1970's by many authors. Research on location problems within the scope of queuing systems was initiated with Larson's paper \cite{larson:74} followed by work of Larson,  Berman and coauthors (e.g. \cite{berman;larson;chiu:85} and \cite{berman;larson;parkan:87}) on discrete location problems. For surveys we refer to the relevant chapters of the collections
\cite{drezner;hamacher:04}[Chapter 11],  
\cite{mirchandani;francis:90}[Chapter 13], (vehicle routing problems under stochastic side constraints),
and the recent very detailed \cite{berman;krass:19}.
We sketch here only two main research directions:\\
(1) In \cite{drezner;schaible;simchi-levi:90} and \cite{scott;jefferson;drezner:99} the authors consider 
location problems where mobile servers move in a plane and  demands of the clients occur as Poisson processes. 
The mobile servers are described as queueing systems,  e.g. of type $M|M|1|\infty$ or  $M|G|s|\infty$.
The guiding principle is to  incorporate travel times to and from the clients into the service times of the mobile server. A survey is provided in \cite{berman;krass:04}.
Locational analysis in a randomly changing environment with occurrence of queueing phenomena is investigated in \cite{dan;marcotte:19}.\\
(2) A class of models which is somehow dual to that described in (1): 
The servers are fixed and clients move to the nearest service stations. Work in this direction is by \cite{berman;drezner:07} and investigated further by \cite{aboolian;berman;drezner:08} and \cite{aboolian;berman;drezner:09}, where additional references are provided.
In these papers the demand  is generated at specified points according to a Poisson or a general renewal process and the servicing nodes usually are modeled as M/M/k/$\infty$ systems. For different cost functions the authors determine the optimal location of a set of servers on a set of nodes.
The field is reviewed in detail in \cite{berman;krass:19} and  classified there as SLCIS ( $=$ Stochastic Location models with Congestion and Immobile Servers).

\underline{Location-inventory problems} (LIPs) aim ``to integrate strategic supply chain decisions with tactical and operational inventory management decisions \cite{farahani;bajgan;fahimnia;kaviani:15}''.
The basic LIP, as described in \cite{farahani;bajgan;fahimnia;kaviani:15}[Section2], encompasses a single supplier (production center),  several warehouses (distribution centers), and retailers.  The location of the supplier and the retailers are given and the decision problem is to determine the location of the distribution centers. \cite{farahani;bajgan;fahimnia;kaviani:15} provides a survey of research on basic LIPs and more evolved variants.

\underline{Location of local repair facilities} with a central service station and inventories at the repair stations is  investigated in \cite{ommeren;bump;sleptchenko:06}.
Given the location of the center, the objective is to determine  optimal locations of repair facilities, ensuring a certain inventory level at each station. The optimization problem is solved by  approximations and local search.
In this article queueing  models are used to describe the interrelation between local repair stations and the center.

\underline{Location routing problems} aim to combine location analysis and planning of vehicle routing.
Fundamentals of vehicle routing problems are described in \cite{laporte:88}.
A recent survey with emphasis on location decisions on networks 
(discrete location problems) is 
\cite{albareda-sambola;rodriguez-pereira:19}. Closer to our problem of location in the plane are \cite{salhi;nagy:09} and in a more versatile setting \cite{manzour-al-ajdad;torabi-salhi:12}.

\underline{Integration of strategic and tactical/operational aspects of planning} is a common topic of almost all the mentioned work. 
\label{page:Integration}
In \cite{salhi;rand:89} for LRPs it is shown that separating decision on location and routing can lead to sub-optimal decisions. 
Under the heading ``Why logistics matters in location modeling'' this problem is discussed indepth in \cite{heckmann;nickel:19}[Section 6.2], stating as main conclusion  ``that making location  decisions ignoring primary logistics activities \dots may result in excessive costs.'' 

\section{Problem description and classification}
\label{sect:Problem}

\underline{Problem setting.} We are given a set of warehouses (stations), indexed by $j\in \{2,3,\dots,J\}$, whose positions
$a_j=(a_{j1},a_{j2})\in\R^2$ in the plane are known. We are also given the aggregate demand (requirement) of $D_j$  truck loads per day generated by warehouse $j$ for a single commodity, $D_j>0$ .\\
We are to find a location $x=(x_1,x_2)\in\R^2$ for a production center (source of commodities), indexed henceforth by $1$, where the commodities are produced and dispatched.\\ 
To deliver the commodities from the production center to the warehouses, $N\geq 1$ identical trucks circulate in the system.
The dispatching rules are roughly as follows:
A truck loaded at the center is directed  to one of the warehouses according to a given schedule or plan, is unloaded there, and returns to the  central station to start another  delivery cycle. The next destination,
 determined by the scheduling regime, might be different.
The center has limitations on its capacity to ship the product, especially when loading several trucks simultaneously, and the warehouses have limited capacity to unload trucks in parallel.\\
To assess the system's performance we take into consideration the distances between the center and the warehouses, the number of trucks, the travel  times,
 the scheduling rules for sending out the trucks,
the loading and unloading times, and  additionally the delay resulting from limitations of loading and unloading capacities, i.e., congestion at loading and unloading facilities.

\underline{The optimization problem} to be solved is to determine the location of the center such that the delivered amount of commodities is maximal and all requests of the warehouses are satisfied with a minimal number of trucks.

\underline{Related problems and classification.} 
Following \cite{min;jayaraman;srivastava:98} and \cite{nagy;salhi:07} this problem is related to location-allocation problems because  it is assumed that there are only radial trips from the center to the warehouses.
The problem is related to location-routing problems because it deals with ``location planning with tour planning aspects taken into account  \cite{nagy;salhi:07}[p.1]'',
because we have to decide about subsequent tours for the trucks to visit different warehouses (the scheduling regimes/plans).
The relation to transportation-location problems in the sense of Cooper \cite{cooper:72}, \cite{cooper:76} is clear because we consider, among others, both problems. Moreover, the center has limitations on its capacity  to ship the product \cite{cooper:72}[p. 94].
Additionally, our optimization criterion encompasses a time dimension as it is considered in \cite{tapiero:71} for the classical setting of Cooper.
Aspects of our problem setting which are not included in the above standard problem classes are:\\
$\bullet$ Limited unloading and loading capacities at the warehouses and the central production facility,  i.e. congestion in the transportation network and queueing problems.\\
$\bullet$ Time dimension of transportation-location-allocation problems, as discussed e.g. in \cite{tapiero:71}[p. 383], in the optimization criterion.
In \cite{tapiero:71}[Section 5] questions concerning delivery time lags are sketched.
In our setting sequentially varied scheduling of trucks is allowed and for feasible scheduling of trucks emerges a necessary condition which optimal plans must satisfy. 
We sketch a toy example which highlights the problems.
\begin{beisp}\label{ex:TimeIncluded}
Consider the case $J=3$, i.e. a center indexed $1$ and two warehouses, numbered  $2, 3$, with demands  
$D_2 = 10$ and $D_3 = 20$ truck loads per day for a commodity and $N=3$ trucks with equal capacity. Two reasonable schedules are\\
\textbf{(i)}
one truck serves station $2$ and two trucks serve station $3$,\\ 
\textbf{(ii)} all trucks are scheduled
to serve in a cycle first station $3$, thereafter station $2$, and finally station $3$ again; these
cycles are iterated.\\
Both schedules generate fair service for the warehouses because they guarantee that over time 
warehouse $2$ will obtain $\rho_2=1/3$  and
warehouse $3$ will obtain $\rho_3=2/3$ of the overall delivered goods.\\
Because we investigate location decisions we have to consider long time horizons and we will therefore assume that the system has approached its stationary state. 
Clearly, then a necessary condition for optimality of a scheduling regime for trucks is to partition in the long run sent-out goods according to: 
$\rho_2=D_2/(D_2+D_3) \wedge \rho_3=D_3/(D_2+D_3)$.
\end{beisp}
Our problem  fits into neither of the  mentioned streams  
on queueing-location problems because we incorporate local congestion at the warehouses explicitly.
Determining a center's location in \cite{berman;drezner:07}, \cite{aboolian;berman;drezner:08}, and \cite{aboolian;berman;drezner:09} is somehow similar to our problem but in these papers no two-way interactions (here generated by trucks) occur between sources of demands  and the production center.

\section{Model and main results}\label{sect:ModelResults}
\subsection{Details of  problem statement}\label{sect:DetailsProblem}

\underline{Modeling the warehouses-production network.}
Warehouse $j$  is equipped with $s_j\geq 1$ service facilities for unloading trucks and  there is  ample waiting space for trucks that arrive while all unloading facilities are busy,
$j\in\{2,\ldots,J\}$. The time for unloading a truck is exponential with rate $\mu_j$ per hour.\\
Center $1$ has $s_1\geq 1$  service facilities for loading trucks and ample waiting space for trucks that arrive while all loading facilities are occupied. Loading times are exponential with mean $\mu_1^{-1}$. 
The queueing regime for trucks is First-Come-First-Served (FCFS).
We abbreviate the service rate functions by
\begin{equation}\label{eq:Muj}
\mu_j(n) := 
\mu_j\cdot min(n, s_j),~~ n=0,1,\dots, \qquad j=1,\dots,J.
\end{equation}
We consider general schemes for scheduling the radial trips of  trucks which must meet only the following restrictions.
If $D:=\sum_{j=2}^J D_j$ denotes the total demand per day, the portion of  demand that has to be delivered to warehouse $j$ is
\begin{equation}\label{eq:routingRho}
\rho_j=D_j/D\in (0,1],~~j\in\{2,\ldots,J\},\quad \sum_{j=2}^J \rho_j = 1.
\end{equation}
Then a truck loaded at the center is directed to warehouse $j$ with (average) frequency $\rho_j=D_j/D$.
In the mathematical model we realize this property of the dispatching rules by a randomized schedule which selects the next warehouse with probability $\rho_j$ for destination $j$
in a Markovian way.
We emphasize that this does not mean that scheduling of trucks should be randomized, but reflects that we are interested in gross characteristics of the system layout. We comment 
on this later in the discussion of ``Random routing\dots'' 
on p. \pageref{page:RandRouting}.\\
Demand that is not satisfied immediately will be backordered at the respective warehouses.
The distance between warehouse $j$ and  center  at $x$
is denoted by $d_j(x):=d(a_j,x)$ where $d:\R^2\times\R^2\longrightarrow\R_+$ is a general convex function.\\
\underline{Details of the optimization problem.}
The aim is to fulfill all demands occurring at the warehouses
with a minimal number of trucks
and to maximize the utilization of the given resources. The latter means that we are to maximize the overall mean number of delivered goods per hour (= time unit for servicing)
which is the sum of the  throughputs of the warehouses  measured in truck loads.
The throughput of a warehouses depends on the coordinates $x=(x_1,x_2)$ of the center and the number of trucks $N$.
Formally:\\
Denote for $T\in\R_+$ by $A_j(N;x)(T)$ the amount of commodities which arrived at  warehouse $j$ and is unloaded
within time horizon $[0,T]$ when the center is located at $x$
and $N$ trucks are cycling.
It will be shown that for $j=2,\dots,J$ the throughput 
\begin{equation}
TH_j(N;x) := \lim\limits_{T\to \infty}\frac{1}{T}A_j(N;x)(T)
\end{equation}
for warehouse $j$ exists 
and determines by  standard ergodicity arguments 
for Markov processes the overall mean number of truck departures per hour from $j$ 
in the stationary system. By stationarity this equals the mean number of truckloads delivered to $j$ per hour.
The total throughput of interest  is 
\[TH_w(N;x)=\sum_{j=2}^{J} TH_j(N;x)\}.
\]
``$ _w$'' indicates that we evaluate only throughputs of warehouses. 
This leads to 
\begin{opt}\label{opt:throughputoptimization} Determine
\begin{eqnarray*}
\min_{N\in\mathbb{N}_+}\left(\max_{x\in\R^2} \Big\{TH_w(N;x)\Big\}\right)
&&~\text{and}~~~~
\arg\Big\langle\min_{N\in\mathbb{N}_+}\left(\max_{x\in\R^2} \Big\{TH_w(N;x)\Big\}\right)\Big\rangle\\
\text{subject to}
&&TH_j(N;x)\geq D_j , j=2,\dots,J.
	\end{eqnarray*}
\end{opt}
It will turn out that the main effort in solving  Optimization Problem
\ref{opt:throughputoptimization} is to solve a sequence of 
maximization sub-problems. These are
\begin{opt}\label{opt:ThroughputMaximization} Determine
	for each $N\geq 1$
	\begin{eqnarray*}
		\max_{x\in\R^2} \Big\{TH_w(N;x)\Big\}
		~~~~\text{and}~~~~
		\arg\Big\langle\max_{x\in\R^2} \Big\{TH_w(N;x)\Big\}\Big\rangle.
	\end{eqnarray*}
\end{opt}
\begin{rem}\label{rem:SmallCapacities}
If loading and unloading capacities are small compared to the demand that has to be delivered, there might be no solution of the Optimization Problem  \ref{opt:throughputoptimization}.
This is due to bottlenecks in the network.
Nevertheless, all Optimization Sub-Problems
\ref{opt:ThroughputMaximization} have solutions.
If necessary, we assume  that the available capacities guarantee that a feasible solution of the problem  exists. 
\end{rem}

\subsection{Main results and detailed analysis}  \label{sect:mainresults}
Facility location in connection  with queueing 
problems usually leads to complex algorithms, see
  \cite{berman;larson;chiu:85}. 
In view of this, our first theorem is counter-intuitive.
\begin{theorem}\label{thm:loesung}
Consider  locations  $a_j=(a_{j1},a_{j2})\in\R^2, j=2,\dots,J,$  in the plane and associated weights $\rho_j$ from \eqref{eq:routingRho}.
    Let $x^{\ast}\in\R^2$ be a solution of the standard Weber problem with weighted distances:
  \begin{equation}\label{weberproblem}
\text{Find}~~    \min_{x\in\R^2}\ \Big\{\sum_{j=2}^J\rho_jd_j(x)\Big\}
~~\text{and}~~ x^{\ast}=\arg\Big\langle\min_{x\in\R^2} \Big\{\sum_{j=2}^J\rho_jd_j(x)\Big\}\Big\rangle.~~~
\end{equation}
    Then $x^{\ast}$
    is a solution of the  Optimization Problem \ref{opt:ThroughputMaximization}
     for any $N\geq 1$  as well.
\end{theorem}
The proof is postponed to Appendix \ref{sect:proofs}.
It relies on the observation that the model for the logistic and services network from Section \ref{sect:Problem}  can be described  in terms of a closed queueing network of Gordon-Newell type. 
The proof of the next theorem will be given implicitly by proving correctness of Algorithm \ref{buzen2} below.
\begin{theorem}\label{thm:MiniMax}
If a solution of the Optimization Problem \ref{opt:throughputoptimization} exists for the capacities $\mu_{j}(\cdot), j=1,\dots,J$, it is uniquely determined if $x^*$ is given.
\end{theorem} 
\begin{rem}\label{rem:StrikingResults}
The results of  Theorem \ref{thm:loesung} and Theorem \ref{thm:MiniMax} are striking, so comments are necessary. \textbf{(i)}
If the 	Optimization Problem \ref{opt:throughputoptimization} 
has a solution, i.e. the side constraints are satisfied with 
capacities $\mu_{j}(\cdot), j=1,\dots,J$, 
these service capacities $\mu_{j}(\cdot)$ (respectively the number of service channels $s_j$) at  warehouses and  center do not matter for optimizing  the overall warehouse throughput  with respect to the location of the center. Similarly, the absolute demands
$D_j$ and the number $N$ of trucks are not relevant for the 
optimal location $x^*$.
 The relevant information for the location decision only comprises \\
$\bullet$ the distances  $d_j(x):=d(a_j,x)$, which determine  travel times, and\\
$\bullet$ the proportions $\rho_j=D_j/D$ of goods to be dispatched to warehouse $j$.\\
\textbf{(ii)}
It is intuitive that increasing the loading capacity at the center  increases throughput at any warehouse.
However, less intuitive is: If we fix the capacities at the center and at all but one dedicated warehouse and increase the unloading capacity at the dedicated warehouse then the throughput at \textbf{all} warehouses increases. Both facts are consequences of Theorem 14.B.13 of
\cite{shaked;shanthikumar:94}.\\
When capacities of loading/unloading facilities change
Theorem \ref{thm:loesung} guarantees that  the decision for the optimal location remains optimal as long as the solution of the Optimization Problem \ref{opt:throughputoptimization} exists.\\
\textbf{(iii)} 
Theorem \ref{thm:loesung} does not propose that warehouse throughput is independent of local properties of the warehouses. Details about functional dependencies will be provided below. Moreover, it is not clear in advance whether a prescribed overall throughput can be met with a given set of parameter values.
If the throughput can be met with the given capacities, Theorem \ref{thm:MiniMax} in connection with the main result of \cite{vanderwal:89} guarantees that by successively adding trucks we can increase the throughput until the total requirements can be dispatched.
Otherwise, if loading/unloading capacities do not suffice,
bottlenecks occur.
Our proofs will show that 
we can increase the throughput by increasing the loading and unloading capacities at the nodes, see Section \ref{sect:NumberOfTrucks}. Theorem \ref{thm:loesung} states that in any case \underline{the selected location remains optimal}.
\end{rem}

\subsection{Analysis of the model as a Gordon-Newell network}\label{sect:reliable stations}
The locations in the problem setting of Section \ref{sect:Problem} can be arranged as a star-like graph with warehouses as exterior vertices $2,\dots,J$ and the production unit as  central vertex $1$. Routes from the center to the exterior nodes, and vice versa back, correspond to edges (links, lanes).
The vertices contain the loading and unloading facilities
modeled as queueing systems.
The circulating trucks are modeled as customers requesting for service at these queueing systems.\\
Because the number of trucks circulating in the network is fixed these features establish a closed queueing network structure (Gordon-Newell network). Additionally,
we apply a standard feature to incorporate travel times into the model for the logistic network:
For each warehouse, roads from the center to that warehouse and back are modeled as two additional infinite server nodes  with random or deterministic service times ( $=$ travel times). 
Necessary  definitions, facts, and formulas from  network theory are summarized in  Appendix \ref{sect:fundamentals}.

\subsubsection{The detailed model of the logistic and services network}\label{sect:DetailedModel}
We start with node $1$ (the center)
and  nodes $2,\ldots,J$ (the warehouses) which are  multi-server nodes with $s_j$ service channels and  exponentially distributed service times with mean $\mu_j^{-1}, j=1, 2,\ldots,J$,
and  with distances $d_j(x):=d(a_j,x), j=2,\ldots,J$.
$N\in\N_+$ customers (trucks)  cycle in the network. For simplicity of presentation we assume that trucks are traveling with unit speed, i.e., $1$ km/hour.
(In examples we shall introduce realistic speeds.)
Whenever a truck is served (loaded) at center $1$ and is routed to warehouse $j$ it has to travel distance $d_j(x)$ to and from.
Traveling these distances is modeled as trucks being served by an infinite server station
for a deterministic time $d_j(x)$ or a random time with mean $d_j(x)$. The first infinite server, from $1$ to $j$, is denoted $ja$. After passing the road to $j$ the truck  will be unloaded at warehouse $j$.
When unloaded  at warehouse $j$ the truck 
travels back to center $1$. This  is modeled as being served by another infinite server station, denoted $jb$, with mean service time $d_j(x)$.\\
Summarizing: Any radial tour from the center to node $j$ consist of three nodes,  representing: (i) traveling to $j$ (via $ja$), (ii) unloading at $j$,
and (iii) return from $j$ (via $jb$). Thus we  have a network with $3 (J-1) + 1$  nodes and  any round trip  (radial tour to $j$) is of the  form: ''$\mathbf{1 \to ja \to j \to jb \to 1}$''.\\
We refer henceforth to center,  warehouses, and lanes jointly as ``stations''.
We assume that all service and travel times are independent and for all lanes and loading/unloading stations identically distributed. The dispatching rules and the radial structure of the trucks are modeled by routing probabilities $r(\cdot,\cdot)$ as follows:
\begin{description}\label{routing}
  \item[(R1)]from center $1$ to lane $ja$: $r(1,ja):= \rho_j$, $j=2,\ldots,J$,~~ $\sum_{j=2}^J r(1,ja)=1$,
  \item[(R2)]all other routing is deterministic: $r(ja,j)=r(j,jb)=r(jb,1)= 1, j=2,\ldots,J$, and
  $r(k,m)=0$ otherwise\,.
\end{description}
Routing decisions at node 1 are independent of the network's previous history.
For simplicity of presentation we assume that the travel times of the trucks are exponentially distributed with mean $d_j(x)$ when $x$ is the location of the center. We  discuss this  in detail in Section \ref{sect:DiscussionModel}.\\
With these stations, routing, and customers we have  constructed a star-like Gordon-Newell network.
Local states of nodes are: For station $j$ is $n_j$ the number of customers present (in service + waiting), for station $ja$ is $n_{ja}$ the number of trucks on the way from station $1$ to station $j$ and
for station $jb$ is $n_{jb}$ the number of trucks on the way from station $j$ to station $1$.
We abbreviate state-dependent service intensities as $\mu_j(\cdot),j=1, 2,\ldots,J,$ given in  \eqref{eq:Muj} and for lanes by
\begin{equation}\label{eq:intensities}
\mu_{ja}(x)(n_{ja}) = d_j(x)^{-1}\cdot n_{ja}~~\text{and} ~~
\mu_{jb}(x)(n_{jb}) =d_j(x)^{-1}\cdot n_{jb},\quad j= 2,\ldots,J.
\end{equation}

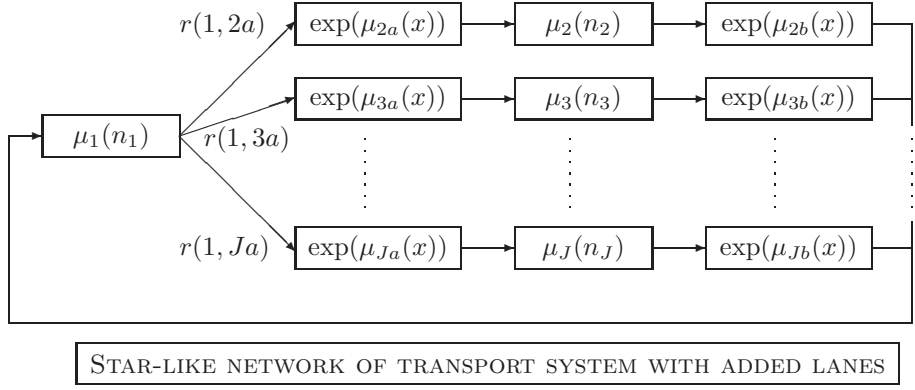
\setlength{\unitlength}{0.9cm}
\begin{picture}(15,5.3)(0.8,-0.4)
    \put(1,2.85){\framebox(2,0.6){$\mu_1(n_1)$}}

    \put(3,4.68){$r(1,2a)$}
    \put(3.35,3){$r(1,3a)$}
    \put(3,1.4){$r(1,Ja)$}

    \put(4.7,4.5){\framebox(2.4,0.6){$\exp(\mu_{2a}(x))$}}
    \put(4.7,3.4){\framebox(2.4,0.6){$\exp(\mu_{3a}(x))$}}
    \put(4.7,1.2){\framebox(2.4,0.6){$\exp(\mu_{Ja}(x))$}}

    \put(7.9,4.5){\framebox(2,0.6){$\mu_{2}(n_2)$}}
    \put(7.9,3.4){\framebox(2,0.6){$\mu_{3}(n_3)$}}
    \put(7.9,1.2){\framebox(2,0.6){$\mu_{J}(n_J)$}}

    \put(10.7,4.5){\framebox(2.4,0.6){$\exp(\mu_{2b}(x))$}}
    \put(10.7,3.4){\framebox(2.4,0.6){$\exp(\mu_{3b}(x))$}}
    \put(10.7,1.2){\framebox(2.4,0.6){$\exp(\mu_{Jb}(x))$}}

    \multiput(7.1,4.8)(2.8,0){2}{\vector(1,0){0.8}}
    \multiput(7.1,3.7)(2.8,0){2}{\vector(1,0){0.8}}
    \multiput(7.1,1.5)(2.8,0){2}{\vector(1,0){0.8}}

    \put(3,3.15){\vector(1,1){1.7}}
    \put(3,3.15){\vector(1,-1){1.7}}
    \put(3,3.15){\vector(3,1){1.7}}

    \put(13.1,4.8){\line(1,0){0.6}}
    \put(13.1,3.7){\line(1,0){0.6}}
    \put(13.1,1.5){\line(1,0){0.6}}

    \put(13.7,4.8){\line(0,-1){1.5}}
    \qbezier[6](13.7,3.1)(13.7,2.6)(13.7,2.1)
    \put(13.7,1.9){\line(0,-1){1.5}}

    \qbezier[6](5.7,3.1)(5.7,2.6)(5.7,2.1)
    \qbezier[6](8.7,3.1)(8.7,2.6)(8.7,2.1)
    \qbezier[6](11.7,3.1)(11.7,2.6)(11.7,2.1)

    \put(13.7,0.4){\line(-1,0){13.2}}
    \put(0.5,0.4){\line(0,1){2.75}}
    \put(0.5,3.15){\vector(1,0){0.5}}
    
    \put(1.5,-0.5){\framebox(12,0.6){\sc Star-like network of transport system with added lanes}}
\end{picture}

Summarizing, we have  a Gordon-Newell network
with node set\\
$\bar{J}=\{1, 2,\ldots, J, 2a,3a,\ldots,Ja,2b,3b,\ldots,Jb\},$
individual service rates $\mu_{ja}(x) = \mu_{jb}(x) =
d_j(x)^{-1}$, and $\mu_{j}(\cdot),$
routing matrix via {\bf (R1), (R2)} (p.\pageref{routing})
and state space
\begin{equation}\label{statespace1}
S(N,{\bar J}) = \{(n_i:i\in {\bar J})\in \mathbb N^{\bar J}:
\sum_{i\in{\bar J}}n_i = N\}.
\end{equation}

\subsubsection{Utilization of resources: Computing throughputs}\label{sect:Throughputs}
Utilizing facts collected in Appendix \ref{sect:fundamentals}, we are in a position to determine  the overall throughputs at the warehouse stations. These are measures for efficient utilization of the given resources. Proofs are postponed to Appendix \ref{sect:proofs}.
\begin{theorem}\label{thm:normingconstant}
Denote by $G(N,\bar{J};x)$ the normalization constant of the  Gordon-Newell network with node set $\bar J$, if $N$ trucks are cycling and the center is located at
$x\in \mathbb R^2$. Then with $\eta_j$ from the routing {\bf (R1), (R2)} via traffic equation \eqref{trafficeqn} it holds with
$\eta_1= \frac{1}{4}, \eta_j= \frac{1}{4}\cdot \rho_j, j=2,\dots, J,$
    \begin{equation*}
        G(N,\bar{J};x)=\sum_{n=0}^N\left[\left(\sum_{n_1+\ldots
        +n_J=N-n}\prod_{j=1}^J\left(\prod_{k=1}^{n_j}\frac{\eta_j}{\mu_j(k)}\right)\right)
        \frac{2^n}{n!}\left(\sum_{j=2}^J
        \eta_jd_j(x)\right)^n \right].
    \end{equation*}
\begin{equation}\label{othruput1}
\text{The overall throughput of the network is}~~    TH(N;x)=\frac{G(N-1,\bar{J};x)}{G(N,\bar{J};x)}.\quad
\end{equation}
The total throughput at the warehouse stations is
\begin{equation}\label{othruput}
 TH_w(N;x) =  TH(N;x)\sum_{j=2}^J \eta_j = \frac{1}{4}\cdot TH(N;x).
\end{equation}
\end{theorem}
\begin{rem}\label{rem:Interpretation}
The representation of $G(N,\bar{J};x)$ has a remarkable interpretation. It is the same normalization constant as that for a Gordon-Newell network with $N$ customers, $J$ multi-server stations with the same service rates  as given in \eqref{eq:Muj} and an attached
\underline{single infinite server}, which will be indexed by $J+1$, with visit ratio $\eta_{J+1} = 1/2$ and exponentially distributed service time with mean $(\sum_{j=2}^J  \rho_jd_j(x))$.
\end{rem}
\begin{kor}\label{normingconstant-det}
Consider the system of Theorem \ref{thm:normingconstant} with deterministic or general random travel times
 with  means $d_j(x)$.
Then the normalization constant is the same $G(N,\bar{J};x)$ and the relevant throughputs are
\eqref{othruput1} and \eqref{othruput} as well.
\end{kor}

\subsection{Determining the number of trucks} \label{sect:NumberOfTrucks}
We  demonstrate the power of Theorem \ref{thm:normingconstant} by showing how to determine efficiently the minimal number of trucks to fulfill the total demand. We assume that 
the center's location $x$ and capacities $\mu_j(\cdot)$  are fixed and sufficiently high to satisfy demands $D_i$ eventually, i.e. with sufficiently many trucks.\\
Recall that in our development we assumed up to now that trucks travel with unit speed (1 km/hour). This implies that $d_j(x)$ is exactly the time to travel distance $d_j(x)$.
For the present demonstration we allow general speed $S>0$ for  trucks. The mean time for traveling  distance $d_j(x)$ is then
the mean service time $d_j(x)/S$ at the infinite servers $ja$ and 
$jb$, $j=2,\dots,J$.\\
With notation from Definition \ref{defn:GN} and
Theorem \ref{thm:gordonnewell}, we apply Buzen's Algorithm
\ref{alg:Buzen} in a first step to a Gordon-Newell network consisting of stations $1, 2, \dots, J$. Then we apply the  representation of $G(N,\bar{J};x)$ from
Remark \ref{rem:Interpretation}.

\begin{alg}
[Determine minimal number of trucks.]\label{buzen2}
Let $C\in(0,\infty)$ denote the capacity of trucks, i.e., the amount of the commodity, that each truck can carry, and
$S\in(0,\infty)$ the speed of the trucks.\\
\underline{Initialization}: {\sc Store}\\
$ G(0,j) := 1, j=1,\dots,J,~~~ \kappa := 2\cdot \sum_{j=2}^J  \eta_jd_j(x)/S,~~~
D=\sum_{j=2}^J D_j$,\\ $H(1,\bar J;x):= 1$.

\noindent
{\sc Set} $N\leftarrow 1$.

\underline{Iterate (*)}
{\sc For} $N$ {\sc do}
\begin{eqnarray*}
   &&\text{{\sc Store}}~~ G(N,1):= g_1(N)~\text{from}~ \eqref{eq:Buzen2}.\\
   &&\text{{\sc Compute with}}~ \eta_1={1}/{4}, \eta_j={1}/{4} \rho_j, ~{\sc and}~
   \mu_j(k) ~{\sc from}~ \eqref{eq:Muj}, j=2,\dots,J,\\
   &&{\sc from}~ \eqref{buzen1}:\quad G(0,J), G(1,J),\dots G(N,J).\\
   &&\text{{\sc Compute}}~ G(N,\bar{J};x)=\sum_{n=0}^N\left[G(N-n,J) \cdot \frac{\kappa^n}{n!}\right].\\
   &&  \text{{\sc If} }~~
C\cdot\frac{1}{4} \frac{H(N,\bar J;x)}{G(N,\bar{J};x)} \geq D\,:\quad \text{{\sc Then output}~~$N$,}\quad \text{{\sc Stop}. }\\
   && \text{\sc{else Store}}\\
   &&{H(N+1,\bar J;x)}\leftarrow{G(N,\bar{J};x)}~
   \text{{\sc and set}}~ N\leftarrow N+1.~~~
   \text{{\sc Go to~~\underline{(*)}}}\\
   &&\text{{\sc Output} is the  minimal number of trucks needed to
guarantee the}\\
 &&\qquad\qquad  \text{required transport capacity to satisfy demand} ~D =D_2+\dots+D_J.
\end{eqnarray*}
\end{alg}
\begin{beweis}
Following Remark \ref{rem:Interpretation},  ${G(N,\bar{J};x)}$ from  Theorem \ref{thm:normingconstant}
can be interpreted as  normalization constant
in a Gordon-Newell network with nodes $1,\dots,J$ having visit ratios $\eta_j$ and
service rates $\mu_j(n_j)=\mu_j\cdot \min(n_j,s_j)$ at station $j$, and an additional infinite server node $J+1$ with visiting ratio $1/2$ and exponential-$(\sum_{j=2}^J  \rho_jd_j(x))$ service time distribution.
Because the $\mu_j(n_j)$ and $n\cdot (\sum_{j=2}^J  \rho_jd_j(x))$ are non-decreasing in $n_j$ and $n$, from van der Wal's theorem \cite{vanderwal:89} it  follows that
the throughput of this artificial Gordon-Newell network is non-decreasing in $N$.
As can be seen from the proof in \cite{vanderwal:89} the throughput is strictly increasing in $N$. This guarantees that the algorithm stops after a finite number of iterations, 
because we assumed that the capacities are high enough to satisfy all the demands eventually.
\end{beweis}
Recall that from Theorem \ref{thm:loesung} the optimal location for the center is independent of loading and unloading capacities and the number of trucks cycling. So, if for   $N$ trucks and optimal location $x$ the  demand $D$ exceeds the  achievable maximal
throughput $TH_w(N;x)$ for the given $\mu_{j}$,  one may increase the loading
and/or  unloading capacities (for a proof see Theorem 14.B.13 of \cite{shaked;shanthikumar:94}).

\begin{rem}
The fact that for a given set of parameters 
$\mu_j(\cdot), j=1,\dots,J,$	it may be impossible to realize the requested demands is a consequence of the observation that in closed queueing networks with nodes which are, roughly stated, not all of infinite server type, bottlenecks exist.
The throughput of bottlenecks converges under unbounded increasing population size to a finite value. This bounds the network's overall throughput which can be attained by increasing number of customers (trucks). 
A short survey with more relevant details of classical bottleneck analysis  is \cite{schweitzer;serazzi;broglia:93}.
\end{rem}

Because the maximal throughput at a bottleneck node can be increased by increasing the local service rate, 
we have the following simple recipe.
\begin{prop}\label{mu1opt}
	Consider the system of Theorem \ref{thm:normingconstant} with all parameters other than $\mu_j, j=1,\dots,J$ being fixed.
	Then loading and unloading capacities $\mu_j$ for  center and  warehouses  exist  which guarantee that
	 $TH_w(N;x) \geq D$ holds for sufficiently many trucks available.
\end{prop}

\subsection{Discussion of the modeling assumptions}\label{sect:DiscussionModel}
\underline{Infinite server queue for modeling traffic on a lane} \cite{newell:82}[Chapter 6] is a  standard device. For fixed mean travel time $d_j(x) <\infty$ of a vehicle moving with unit speed, we can allow any distribution to incorporate or forbid overtaking. 
Extreme travel time distributions are  exponential-$d_j(x)^{-1}$ (maximal entropy) and deterministic-$d_j(x)$ (minimal entropy). In any case:
The joint stationary queue length distribution of number of trucks on lanes, normalization constants, and  throughputs  remain the same.
A realistic model for travel times is obtained using random  travel times with mean $d_j(x)$ and small variance which generates moderate overtaking.
We obtain  throughput  \eqref{othruput1} and \eqref{othruput} in any case.
We discuss consequences of this observation  in Section \ref{sect:Sensitivity}\\
\underline{Exponential  multi-server stations under FCFS}  for loading and unloading the trucks are standard models. Other service disciplines may be more realistic in specific situations.
These  often yield the same performance characteristics with respect to the optimization
criterion, see \cite{daduna:01a}[Theorem 9.9].
Moreover, for many service disciplines it is possible that loading and unloading times may have general distributions.
Two typical settings where throughput will be the same for any shape of the service time distributions as long as the mean service time is fixed (robustness of throughput), are (see \cite{daduna:01a}[Theorem 9.7 (3),Theorem 10.2, Remark 9.6]):\\\label{page:ServiceDisz}
(i) If there is ample  capacity, i.e. all trucks present are served  in parallel (= infinite server),  general distributions for loading and unloading times are admitted.\\
(ii) If loading of all trucks present at the central station is is performed concurrently, the adequate model of servicing  is Processor Sharing. This means: If there are $n_1$ trucks at the center each of them obtains a fraction of $1/{n_1}$ of the station's total capacity.\\
An interesting observation is that the  mentioned robustness  property (insensitivity) is not valid in the related
{\em combined location-routing problem} (LRP) where strategic and tactical decisions are intertwined. For more details 
see the review paper \cite{nagy;salhi:07}[Section 1.2] and the early survey \cite{laporte:88}.\\
\underline{Location of the center coincides with one of the warehouse locations}  is a possible scenario
because we have reduced the queueing-location problem to a pure location problem. If the center's position is $x=a_j$, then the travel times to and from this node are $d_j(x)=0$. Consequently, the service times at the infinite servers  $ja$ and $jb$ are zero. This situation is covered by our framework with mean zero travel times from $1$ to $j$ and back. In this case $n_{ja} = n_{jb} =0$.\\
\underline{Random routing when departing from the center} 
\label{page:RandRouting} is a modeling assumption which 
adjusts the distribution of the available transportation capacity in the long run and in the stationary system according to the demands of the warehouses. As discussed in Example \ref{ex:TimeIncluded}, transitions into realizable schedules can be found easily.
A realistic schedule  determines a sequence of visits for each truck to warehouses. These sequences are cyclically iterated.
A necessary condition for optimal scheduling
is that (in a stationary system and in the long run) the  schedule  \underline{must} reproduce the
frequencies $\rho_j=D_j/D$:\\
Recall Example \ref{ex:TimeIncluded}: Two warehouses with demands $D_2 = 10$ and $D_3 = 20$ truck loads  per  hour and $N=3$ trucks with equal capacity. 
Both reasonable schedules described there can be modeled with deterministic routing schemes described in
\cite{kelly:79}[Section 3.4] to obtain a Markovian network model using a more elaborated state space.
The route of the trucks (in the sense of Kelly) for the second schedule would be:\\
 $[1\to 3a\to 3\to 3b\to 1\to 2a\to 2\to 2b\to 1\to 3a\to 3\to 3b] \to [1\to 3a\to3b\to 1\ldots\to 3b]\to[\ldots]\to \dots$,
iterated indefinitely over time.\\
If we evaluate the joint queue length distribution of stations $1, 2, 3,$ in these models we  obtain
exactly the stationary distribution of our present model with $\rho_2 = 1/3, \rho_3 = 2/3$.\\
The main conclusion from  Theorem \ref{thm:loesung} is:
The optimal location of the center is the same \textbf{for all schedules which generate in the stationary state the values $\rho_j=D_j/D$ by cyclical sequencing}. This is in line with intuition.\\
\underline{Incorporating time into the optimality conditions,}
 we follow the arguments of Tapiero \cite{tapiero:71}:
``The time dimension in problems of transportation-location-allocation is particularly important since decisions to construct production facilities are based on long range plans. 
Also, although environmental conditions, demand, etc. may change over time, the decision to locate a plant in a particular place is made once and is not subject to frequent change. (page 383)''\\
``Long range plans'' in the sense of Tapiero \cite{tapiero:71} justify to assume that we consider a
\textbf{stationary} system. Note that this does not 
mean to consider static systems. Only the random fluctuations of the system are invariant over time.\\

\subsection{Minimizing round-trip times} \label{sect:MiniPassageTimes}
Throughput maximization is probably the most important objective
in the logistic and services network. Another important objective is minimization of expected round-trip times, i.e.
the mean  travel time for a truck between two successive departures from the center.

The \it expected passage-time $Z_j(N;x)$ from the center located at $x$ to station $j$ and back \rm is the sum of all expected waiting and service times, which a truck spends
at station $j$  plus traveling to and from $j$
and thereafter at station $1$. If $W_i(N;x)$ is the expected sojourn time at station $i$, then
$Z_j(N;x)=W_{ja}(N;x)+W_{j}(N;x)+W_{jb}(N;x)+W_{1}(N;x),\ \ j=2,\ldots,J.$
The overall \it expected passage-time \rm is
   $ Z(N;x):=\sum_{j=2}^J r(1,ja)Z_j(N;x).$\\
With $\eta_j=(1/4)r(1,ja)$, $j=1,\ldots,J$ and
$\eta_1=\sum_{j=2}^J\eta_j$ (see \eqref{trafficsolution}) we get
\begin{equation}\label{ZNX1}
    Z(N;x)=\sum_{j=2}^Jr(1,ja)Z_j(N;x)
    =4\sum_{j=2}^J\eta_jZ_j(N;x)
    =4\sum_{j\in\bar{J}}\eta_jW_j(N;x).
\end{equation}
From  Little's Theorem  \cite{chen;yao:01}[Formula (2.18)], we obtain
    $\sum_{j\in\bar{J}}\eta_jW_j(N;x)=\frac{N}{TH(N;x)}.$
This yields  the optimization problem: 
\begin{displaymath}
\text{Find}\quad    \arg\min_{x\in\R^2}\ \left\{Z(N;x)=4\frac{N}{TH(N;x)}\right\}.
\end{displaymath}
So $Z(N;x)$ attains its minimum when $TH(N;x)$ is
maximal and  $Z(N;x)$ attains
its minimum at $x^\ast$ given in Theorem \ref{thm:loesung}. Hence,  travel time minimization is reduced to a standard Weber problem.
This result holds for the generalizations from Section \ref{sect:DiscussionModel} as well because we consider  mean passage times.

\section{Numerical examples and discussion} \label{sect:NumericManagerSensitive}
Our main result states: For the Weber problem in  the logistics and services network under congestion the strategic decision for the center's location and the 
tactical/operational decision for the fleet size decouple
as long as the relative demands $\rho_j:=D_j/D$ remain stable.
Nevertheless, it is of value to demonstrate the consequences of this invariance property by  examples and to discuss  consequences of the result for managerial decision making.
Distance measure is in any case Euclidean distance  with specific weights which will vary.

\subsection{Numerical example}\label{sect:Numeric}
The first two examples mimic the location of $12$ midsize up to large towns in Northern Germany in a rectangle of size 400 km $\times$ 260 km (approximately).
We embed this rectangle into the positive lattice
$\mathbb{N}_0\times\mathbb{N}_0$ in the plane, and shift the town in the south-west corner to the point $(10,10)$.
The demands $D^{(\cdot)}_j$ (measured in truck loads per day) of the warehouses $j=2,\dots,J$, are
approximately chosen 
(i) proportional to the number of inhabitants  of the respective cities, $(D_j^{pro})$,
resulting in total demand of $81$ truck loads/day, and
(ii) according to the logarithm of the number  of inhabitants (divided by 1000) of the respective cities, $(D_j^{log})$, resulting in total demand of $66$ truck loads/day.\\
We assume for simplicity of computations that all loading and unloading facilities are single servers.
Unload capacities at all warehouses are  
$2$ truck loads/hour ($\mu_j=2$), loading capacity at the center is $4$ truck loads/hour ($\mu_1=4$).
The locations of the warehouses $(a_{j1},a_{j2}),j=2,\dots,J,$ are indicated in the next table. Demands of the respective locations are listed below the locations. (All numbers are rounded to integers).\\
\noindent
\small{
\begin{tabular}{|r|*{12}{r|}} \hline
 $j =$& 2& 3& 4& 5& 6& 7& 8& 9& 10& 11& 12& 13\\ \hline	
$a_{j1}$ & 10& 100& 170& 290& 410& 220& 260& 180& 320& 160& 40& 80\\ \hline
$a_{j2}$ & 10& 130& 190& 30& 70& 230& 190& 270& 250& 50& 40& 180\\ \hline
$D_j^{pro}$& 3& 6& 19& 2& 36& 2& 1& 2& 2& 5& 2& 1\\ \hline
$D_j^{log}$& 6& 6& 7& 5& 8& 5& 4& 5& 5& 6& 5& 4\\ \hline  
\end{tabular}
}
\\

We applied the Weiszfeld algorithm 
and obtained 
(i) for weights $D_j^{pro}$ the location for the center 
at $x=(x_1,x_2) = (288.156, 112.283)$,
and (ii) for weights $D_j^{log}$ the location for the center 
at $x=(x_1,x_2) = (179.756, 155.904)$.
Without weights the center's location is $x=(179.210, 162.372)$.
We observed the following results of the optimization procedure.\\
(i) For demands $D_j^{pro}$ the distance between the center  with weights and the center without weights is 
$119.909$. The necessary number of trucks needed to fulfill the total demand of $81$ truck loads is $28$ for center with weights and $29$ for center without weights.
Nevertheless, the thoughput at the warehouses is greater with weights and 28 trucks than without weights and 29 trucks. The details are summarized in the following table. For completeness we added the probability that the loading server at the center is busy as a measure of congestion.\\

\begin{tabular}{|r|r|r|r|r|r|}\hline
	demand $D_j^{pro}$&location,~ $\mu_1=4$&trucks&throughput/day& $P(X_1>0)$ \\ \hline
	with weights& 288.156, 112.283 & 28 & 82.261&                0.857\\ \hline
	no weights& 179.210, 162.372 & 29&81.342&                    0.847\\ \hline
\end{tabular}
\vspace{0.3cm}\\
(ii) For demands $D_j^{log}$ the distance between the center with weights and the center  without weigts is 
$6.491$ and the necessary number of trucks needed to fulfill the total demand of $66$ truck loads is the same when the center's location is determined with or without weights. Nevertheless, the throughput at the warehouses is slightly greater with weights. Details are given in the next table.\\ 

\begin{tabular}{|r|r|r|r|r|r|}\hline
demand $D_j^{log}$&location,~$\mu_1=4$&trucks&throughput/day& $P(X_1>0)$ \\ \hline
with weights& 179.756, 155.904 & 19 & 67.871&   0.706990             \\ \hline
no weights& 179.210, 162.372 & 19&67.841& 0.706676                   \\ \hline
\end{tabular}
\vspace{0.3cm}
\\
Interpretation: 
The demand structure $D_j^{log}$ is rather homogeneous and therefore the difference between the locations of the center is  insignificant.
The  case of demand structure $D_j^{pro}$ is more interesting because of the great distance between the centers' location
which is a consequence of the more variable demand structure.
In any case the congestion, measured as $P(X_1>0)$, increased with increasing throughput which follows from easy computations. Our experiments show that the increase of congestion seems to be tolerable, especially if we can reduce the number of trucks in parallel.\\
It is easy to see that in both scenarios the loading server at the center is the bottleneck of the network.
We therefore investigated the influence of the loading capacity at the center in the above examples
and reduced the loading capacity from $\mu_1=4$  to $\mu_1=3$ truck loads/hour. The results for demands $D_j^{log}$ in the
next table are not surprising. The reduced server capacity is compensated by more trucks to deliver total demand of 66.\\

\vspace{0.02cm}
\begin{tabular}{|r|r|r|r|r|r|}\hline
	demand $D_j^{log}$&location,~$\mu_1=3$&trucks&throughput/day& $P(X_1>0)$ \\ \hline
	with weights& 179.756, 155.904 & 22 & 67.054&   0.931308          \\ \hline
	no weights& 179.210, 162.372 & 22 &67.040& 0.931110                   \\ \hline
\end{tabular}
\vspace{0.3cm}\\
For the more variable demand $D_j^{pro}$  (total 81) the results are given in the next table. Because the center is the bottleneck the results for throughput (should be  less than $72$) and non-idling probabilities (should be less than $1$) are due to rounding errors.\\

\vspace{0.02cm}
\begin{tabular}{|r|r|r|r|r|r|}\hline
	demand $D_j^{pro}$&location,~$\mu_1=3$&trucks&throughput/day& $P(X_1>0)$ \\ \hline
	with weights& 288.156, 112.283 & -- & 72.000&                1.000\\ \hline
	no weights& 179.210, 162.372 & --&72.000&                    1.000\\ \hline
\end{tabular}
\vspace{0.3cm}\\
The algorithm stopped when it detected that $100$ trucks are not sufficient to satisfy the requested total demand of $81$ truckloads.
This is a consequence of the fact that the departure stream from the bottleneck server approaches asymptotically   a Poisson process with intensity $3/hour$ which determines asymptotically (Number of trucks $\to \infty$) the maximal  total throughput of the system. The resulting upper bound for the total throughput ($72$ truckloads) does not meet the requested total demand of $81$.
Increasing the service rate at the center to $3.38/hour$ yields
the following results.\\

\begin{tabular}{|r|r|r|r|r|r|}\hline
	demand $D_j^{pro}$&location,~$\mu_1=3.38$&trucks&throughput/day& $P(X_1>0)$ \\ \hline
	with weights& 288.156, 112.283 & 43 & 81.013&                0.998676\\ \hline
	no weights& 179.210, 162.372 & 45& 81.021&                    0.998780\\ \hline
\end{tabular}
\vspace{0.3cm}\\
Note, that overshot of throughput is slightly higher without weights and $45$ trucks  than in case of weights with $43$ trucks needed to satisfy demand.
The moderate deviation of the number of trucks needed in the setting with and without weights is surprising.
But from observations in a series of experiments we concluded that this is not unusual.

Additionally, we performed experiments with $12$ warehouses located on the lattice 
$\{10,\dots,410\}\times\{10, \dots,270\}$  and sampled  independently according to
uniform distribution. The respective demands are selected according to uniform distribution $U(A)$ on different sets $A$ of feasible demands. In any case the unloading capacities are  $\mu_j=2~\text{truck loads/hour}, j= 2,\dots,13$.
The loading capacities $\mu_1$  varied and are given in 
Table \ref{tab:Summary} below.  Details are presented  in the tables of Section \ref{sect:AddNumeric}.\\
For comparison with the ratio (sample mean/sample variance) of 
$D_j^{log}, j=2,\dots,13: (5.5/1.37)$ and	$D_j^{pro},  j=2,\dots,13 : (6.75/109.30)$ in the previous experiments
we indicate for the respective  demand distribution the mean (Exp) and variance (Var) in the first two rows of Table \ref{tab:Summary}.\\
In any of 4 blocks ((I),\dots,(IV)) we performed 10 experiments.
For extreme total demands  it turned out that the demand can not be satisfied with the given capacities due to occurrence of bottlenecks (the number of such samples are indicated as ``Num$\infty$Tru''). We included these cases in the tables in Section \ref{sect:AddNumeric}  and performed additional experiments to obtain in any block 10 complete data sets. 
Within these  we observed throughout that in approximately half of the samples  for the center's location with weights less trucks (usually  1 truck less) are needed than without weights (precise numbers indicated as ``Num$<$Tru'').
We observed in both blocks (III) and (IV) a single experiment where 3 trucks less are sufficient (underlined). These high differences coincided with the maximal distance 
(MaxDist) between the two centers' location (bold).\\
In all experiments where  an equal number of trucks is needed for center with weights and for center without weights, this coincides with a higher throughput at the warehouses for the center's location selected with weights. This seems to be in line with intuition.
On the other side, when less trucks are needed for center with weights, in almost all cases the throughput with more trucks and location selected without weights produces more throughput, i.e. more overshot.\\
This demonstrates that  decisions on the basis of integrated models leads to better utilization of the given resources. This is substantiated by the following observations:\\ 
Exceptions of the moderate decrease of needed resources (trucks) are the extreme cases in (III) and (IV). With 3 trucks less the throughput with weights exceeds the throughput without weights. Moreover, in both blocks 
(III) and (IV) the second largest distance (bold) between the centers' location generated with 1 truck less (with weights)  a higher throughput than without weights.
In Table \ref{tab:Summary} 
we report for any block the minimal and maximal distance between the centers and the minimal and maximal demands.
Detailed results are presented in Section \ref{sect:AddNumeric} of the Appendix.\\
\begin{table} 
	\begin{tabular}{|r|r|r|r|r|}\hline
		&(I)&(II)&(III)&IV\\ \hline
		&$U(\{1,..,8\})$
		&$U(\{1,..,16\})$ & $U(\{1,..,21\})$ & $U(\{1,11, 21\})$\\
		\hline
		Exp& 4.75 & 8.5& 11& 11  \\ \hline
		Var& 5.25& 21.25& 36.67& 66.67  \\ \hline
		$\mu_1$&  4&  5&  7&  7  \\ \hline
		MinDist& 7.5896& 16.2628 & 11.0563 & 5.5115	
		\\ \hline
		MaxDist&60.7717	& 88.8410 & 108.5533 & 85.6094
		\\ \hline
		MinDem& 41 & 93  & 94 &	42
		\\ \hline
		MaxDem&  63  & 123  & 163&  162	
		\\ \hline
		Num$<$Tru& 4 & 6 & 5 & 6	
		\\ \hline
		Num$\infty$Tru& 0 & 2 & 1 & 2	
		\\ \hline
	\end{tabular}
\caption[Summary]{Summary of the experiments. Details in Section \ref{sect:AddNumeric}}. \label{tab:Summary}
\end{table}

\subsection{Sensitivity analysis and robustness}\label{sect:Sensitivity}
Our main results indicated that some of the 
system's parameters  are not relevant for the decision problems which are in the focus of our investigations.

\textbf{(a)} Discussing modeling assumptions in Section
\ref{sect:DiscussionModel} we indicated that
for fixed mean travel time $d_j(x) (<\infty)$ of  vehicles we can allow any shape of travel time distribution.
The joint stationary queue length distribution, the normalization constants, and the throughputs of the
network are the same.\\
Consequently, our model is robust against changes of these data, e.g. against variability of travel times. 
This flexibility is due to so-called  {\it Insensitivity theory for queueing networks}.
This theory dates back to insensitivity in
{\it Verallgemeinerte Bedienungsschemata} \cite{koenig;matthes;nawrotzki:74} and in BCMP and Kelly networks \cite{schassberger:78a}, for details 
see \cite{daduna:01a}[Section 9 and 10].
The relevant fact from insensitivity theory for our problem is:
{\it At an {\em infinite server} the stationary queue length distribution is {\em invariant} under variation of the 
shape of the service time distribution as long \underline{as the  mean is fixed}.}
This implies that we can compute  stationary queue length distributions on the lanes using exponential service time - the result is the same for any other distribution with the same mean, see \cite{daduna:01a}[Theorem 9.7 (3) and Theorem 10.2] 
and the remark thereafter.\\
A similar robustness property is observed when varying service
disciplines (i.e. reorganizing loading/unloading) at the stations to a certain extend (see p. \pageref{page:ServiceDisz}).

\textbf{(b)} Separability of Optimization Problem
\ref{opt:throughputoptimization} implies that decision for the location of the central production facility is robust against variations of all parameters of the integrated production-transportation-inventory system as long as the proportions of the demands $\rho_j=D_j/D, j=2,\dots, J,$ are not changed and the capacities $\mu_j(\cdot)$ are sufficiently high to meat the demand.

\textbf{(c)} On the other side, separability of Optimization Problem \ref{opt:throughputoptimization} implies that, when
the optimal location is fixed,
we can optimize for the number of trucks needed to satisfy  demands $D_2, \dots, D_J$ by reallocation of capacities at the unloading service stations or by adding capacity at the loading station at the center.

\subsection{Managerial insights}\label{sect:ManagerialDecision}

Theorems \ref{thm:loesung} and \ref{thm:MiniMax} substantiate conclusions for managerial decision making.

\textbf{(1)} Strategic (location of central facility) and tactical and operational (routing, scheduling) decisions are usually thought to be independent and are consequently 
separated. Then decision about the location of a central production facility neglects the actual and future capacities, the resulting congestion, and delays at downstream warehouses. Several authors have shown that for location-routing problems (LRPs) such a structural separation produces sub-optimal solutions for allocation problems, see e.g. \cite{nagy;salhi:07}.
For general supply chain analysis this is discussed in \cite{heckmann;nickel:19}.\\
Contrary to this, our Theorems \ref{thm:loesung} and  \ref{thm:MiniMax} justify in a stylized but rather general model the separation of strategic decisions for locations from  several tactical and operational decisions for scheduling and placing service capacities. In our setting, in a first step the  optimal location can be determined by way of a standard Weber problem where future allocation/routing decisions are incorporated only via gross information (or gross assumptions) about expected demand.
Later on, in a second step the fine-tuning of the radial trips can be carried out according to the service resources at hand  without  the central location becoming sub-optimal.
In Section \ref{sect:DiscussionModel} it is demonstrated that results obtained in the stylized model are valid in more realistic settings as well.
  
\textbf{(2)}  Once the system is built, the quality of the delivering process measured in the standard metrics throughput or round-trip time can be increased  by local enhancement of service without making the central location sub-optimal.

\textbf{(3)}  Shifting capacities between the nodes is possible without perturbing  optimality of the center's location as long as the fractions of demand at the stations remain the same. More precisely: If we  can fine-tune
 scheduling at the warehouses by
placing a prescribed  number  $M>J$ of service facilities at the exterior stations, with at least one facility per station, according to some
further optimization criterion, our theorems state that the location of the central server remains optimal as long as $M$ provides a feasible distribution of capacities for loading and unloading.
This second step of fine-tuning is related  to 
distribution of servers in the Multiple Server Location problem introduced in \cite{berman;drezner:07} and investigated further
in \cite{aboolian;berman;drezner:08}, \cite{aboolian;berman;drezner:09}.

\textbf{(4)} Concurrent optimization for location of central facility and the  number of trucks in the system for a target throughput can be carried out in a step-by-step procedure according to Theorem \ref{thm:MiniMax}:
With the center's location fixed  we can algorithmically solve for the needed number of trucks,
see Section \ref{sect:NumberOfTrucks}.\\

\section{Conclusion and directions of further research}\label{sect:Conclusion}
We have developed a methodology for determining jointly
optimal solutions for location-allocation-routing problems which are usually considered to be problems on different levels of decision making: strategic versus tactical/operational level.
We translated the problem into a stochastic network problem and showed that (i) determining the center of a star-like network under constraints on the demands generated by the exterior nodes, and (ii) determining the optimal number of customers in the network, can be separated.\\
Starting from an exponential version of the problem, which allows for simple proofs, we have shown that more realistic models are covered by the result.
Probably, the most important result  is that the decision upon the location of the production center can be decoupled from building the exterior stations and their equipment.\\
Many research problems are not yet tackled in the more general area of locating additional nodes in networks of queues. These can be easily identified by investigating the more involved location problems, like p-median or p-center problems, in the continuous as well as in the discrete setting of queueing networks.\\
Parts of our ongoing research related to the present paper are
(i) location theory in the discrete queueing network setting with prescribed network graphs,
and (ii)  p-center location problems in the setting of this paper.

\begin{appendix}
\section{Appendix}
\subsection{Prerequisites  from  queueing
network theory}\label{sect:fundamentals}
\begin{defi}\label{defn:GN}
     A \textbf{Gordon-Newell network} consists of stations  $\{1,2,\ldots,I\}$. Station j has $s_j\geq 1$ service channels and ample waiting room under FCFS.
    $N>0$ indistinguishable customers
    cycle according to an irreducible Markov matrix
    $R=(r(i,j),i,j=1,\ldots,I)$ in the network and request for service  the nodes. The service time at node $j$ is
    exponentially distributed with mean $\mu_j^{-1}$.
    Whenever $n_j$ customers are present at node $j$ (in service or waiting),  service is provided with rate
   $ \mu_j(n_j) = \mu_j\cdot \min(n_j,s_j)\,.$\\
   Let $X_j(t)$ denote the number of customers  at station j at time $t\ge 0$ and
$X(t):=(X_j(t),j=1,\ldots,I)$  the joint queue
    length vector at time $t$. $X=(X(t),t\ge 0)$ is the joint queue length process on  state space
      $S(N,I)=\{(n_1,\ldots,n_I)\in\N^I, n_1+\ldots +n_I=N\}.$
\end{defi}
\begin{theorem}\label{thm:gordonnewell}
    \bf (\cite{jackson:63}, \cite{gordon;newell:67}) \it
    The joint queue length process $X=(X(t):t\geq 0)$ of the
    Gordon-Newell network is an ergodic Markov process.
    Denote by
    $\eta=(\eta_1,\ldots,\eta_I)$ the unique probability
    solution of the traffic equation
    \begin{align}\label{trafficeqn}
        \eta_j=\sum_{i=1}^{I}\eta_i r(i,j),\ j\in \{1,\dots,I\}.
    \end{align}
     $\eta_j$ is the customers' visit ratio at node $j$.
    With normalization constant $G(N,I)$, the unique stationary and limiting distribution $\pi=\pi(N,I)$
    of X on S(N,I) is
    \begin{align*}
        \pi(n_1,\ldots,n_I)=G^{-1}(N,I)\prod_{j=1}^I\prod_{i=1}^{n_j}\frac{\eta_j}{\mu_j(i)},\quad (n_1,\ldots,n_I)\in S(N,I).
    \end{align*}
\end{theorem}
\textbf{Remark.} Any non-zero solution $\eta$ of \eqref{trafficeqn} is admissible to compute $\pi$. Consequently, the next algorithm can be used with any such $\eta$.
\begin{alg} \label{alg:Buzen}\textbf{Buzen's Algorithm.} 
	\cite{bruell;balbo:80}[Section 2.2.1]
\begin{eqnarray}
	&&\text{For}~j=1,\dots,J~\text{set}~ g_j(0) := 1, ~~ g_j(n_j) := \prod_{k=1}^{n_j}\frac{\eta_j}{\mu_j(k)}, n_j\geq 1.\label{eq:Buzen2}\\
	&&\text{Set boundary values}~~\nonumber\\
	&&   G(0,j) := 1, ~~j= 1,\dots,J,\quad G(m,1) := g_1(m),~~ m=1,2,\dots.\nonumber\\
	&&\text{Denote for}~~j\geq 1, m\geq 1,\nonumber\\
	&& G(m,j) := \sum_{n_1+\dots +n_j=m}\prod_{j=1}^j\left(\prod_{k=1}^{n_j}\frac{\eta_j}{\mu_j(k)}\right)
	=\sum_{n_1+\dots +n_j=m}\prod_{j=1}^j g_j(n_j)\,.\nonumber
\end{eqnarray}
{\em Buzen's algorithm to compute norming constants}
for $j\geq 1$ and $m\geq 1$ is
\begin{equation}\label{buzen1}
G(m,j) = \sum_{\ell=0}^{m} G(\ell,j-1) \cdot g_j(m-\ell)  \,.
\end{equation}
\end{alg}

\begin{lem}
    \bf\cite{chen;yao:01} \it In a Gordon-Newell network with $N\ge 1$ customers the mean number of departures per time unit from node $j$ (node-$j$ throughput) is
    \begin{align}\label{localthru}
        TH_j(N):=\sum_{(n_1,\ldots,n_I)\in S(N,I)}\pi(n_1,\ldots,n_I)
        \mu_j(n_j)=\eta_j&\frac{G(N-1,I)}{G(N,I)}.\\
    \label{globalthru}
\text{The (overall) throughput is}~~  TH(N)=\sum_{j=1}^I TH_j(N)=&\frac{G(N-1,I)}{G(N,I)}.
    \end{align}
\end{lem}

\subsection{Proofs}\label{sect:proofs}
Recalling  node set $\bar{J}=\{1,2a,3a,\ldots,Ja,2,3,\ldots,J,2b,3b,\ldots,Jb\}$, 
the probability solution of the traffic equation \eqref{trafficeqn} for the network in Section \ref{sect:reliable stations} is
\begin{align}\label{trafficsolution}
\eta_1= \frac{1}{4}\ \ \mbox{and \ } \eta_j=\eta_{ja}=\eta_{jb}=\frac{1}{4}r(1,ja)
=\frac{1}{4}\rho_j, \quad    j=2,\ldots,J.
\end{align}
The stationary distribution of the network is
$   \pi(n_j,j\in\bar{J};x) =$
\begin{align*}
    &G^{-1}(N,\bar{J};x)\prod_{j=1}^J
    \left(\prod_{k=1}^{n_j}\frac{\eta_j}{\mu_j(k)}\right)
    \prod_{i=2}^J\left(\prod_{l=1}^{n_{ia}}\left(\frac{\eta_{ia}}{l\mu_{ia}(x)}\right)
    \prod_{m=1}^{n_{ib}}\left(\frac{\eta_{ib}}{m\mu_{ib}(x)}\right)\right)\\
    &=G^{-1}(N,\bar{J};x)\prod_{j=1}^J\left(\prod_{k=1}^{n_j}\frac{\eta_j}{\mu_j(k)}\right)
    \prod_{i=2}^{J}\frac{1}{n_{ia}!n_{ib}!}(\eta_id_i(x))^{n_{ia}+n_{ib}}.
\end{align*}
where we utilized $\eta_{ja}=\eta_{jb}=\eta_j$ and
$\mu_{ja}(x)=\mu_{jb}(x)=d_j^{-1}(x)$ for $j=2,\ldots,J$.
The following representation of normalization constants will be of value.
\begin{lem}\label{lem:normingconstant}
    The normalization constant of the system is
    \begin{align}
        G(N,\bar{J};x)&=\sum_{n=0}^N\left[\left(\sum_{n_1+\ldots
        +n_J=N-n}\prod_{j=1}^J\left(\prod_{k=1}^{n_j}\frac{\eta_j}{\mu_j(k)}\right)\right)
        \frac{2^n}{n!}\left(\sum_{j=2}^J
        \eta_jd_j(x)\right)^n \right].\nonumber\\
   &\text{With}~~~     C_n(N,\bar{J}):=\left(\sum_{n_1+\ldots +n_J=N-n}
        \prod_{j=1}^J\left(\prod_{k=1}^{n_j}\frac{\eta_j}{\mu_j(k)}\right)\right)
        \frac{2^n}{n!}\nonumber\\
   &\text{  and}~~~      h_j(x)  :=\eta_jd_j(x),\qquad
        h(x)  :=  \sum_{j=2}^J h_j(x)=\sum_{j=2}^J\eta_jd_j(x)\nonumber\\
  &\text{we can write}~~~       G(N,\bar{J};x)=\sum_{n=0}^N C_n(N,\bar{J})h(x)^n. \label{ofvalue}
    \end{align}
\end{lem}
\begin{beweis}
   From the definition, we have $G(N,\bar{J};x)=$
\begin{align*}
    &=\sum_{n_1+\ldots +n_{Jb}=N}\prod_{j=1}^J\left(\prod_{k=1}^{n_j}\frac{\eta_j}{\mu_j(k)}\right)\prod_{i=2}^{J}\frac{(\eta_id_i(x))^{n_{ia}+n_{ib}}}{n_{ia}!n_{ib}!}\\
               &=\sum_{n=0}^N\left[\left(\sum_{n_1+\ldots+n_J=N-n}\prod_{j=1}^J \left(\prod_{k=1}^{n_j}\frac{\eta_j}{\mu_j(k)}\right)\right)\left(\sum_{n_{2a}+\ldots +n_{Ja}+\atop n_{2b}+\ldots +n_{Jb}=n}\prod_{j=2}^J\frac{h_j(x)^{n_{ja}+n_{jb}}}{n_{ja}!n_{jb}!}\right)\right].
 \end{align*}
    The statement follows from
    \begin{eqnarray*}
        && \sum_{n_{2a}+\ldots +n_{Ja}+\atop n_{2b}+\ldots +n_{Jb}=n}\prod_{j=2}^J
            \frac{h_j(x)^{n_{ja}+n_{jb}}}{n_{ja}!n_{jb}!}\\
            &=&  \sum_{n_{2a}+\ldots +n_{Ja}+\atop n_{2b}+\ldots +n_{Jb}=n}\prod_{j=2}^Jh_j(x)^{n_{ja}}\prod_{i=2}^Jh_i(x)^{n_{ib}} \frac{1}{\prod_{i=2}^Jn_{ia}\prod_{i=2}^Jn_{ib}}\\
        &=&\sum_{n_{2a}+\ldots +n_{Ja}+\atop n_{2b}+\ldots +n_{Jb}=n}\underbrace{\prod_{j=2}^J\left(\frac{h_j(x)}{2h(x)}\right)^{n_{ja}}\prod_{i=2}^J\left(\frac{h_i(x)}{2h(x)}\right)^{n_{ib}} \frac{n!}{\prod_{i=2}^Jn_{ia}\prod_{i=2}^Jn_{ib}}}_{\mbox{density function of a multinomial distribution}}\cdot\frac{(2h(x))^n}{n!}\\
        &=&\frac{2^n}{n!}h(x)^n.
    \end{eqnarray*}
\end{beweis}
\begin{beweis} (\underline{of Theorem \ref{thm:normingconstant}}) The representation of the normalization constant is the first statement of Lemma \ref{lem:normingconstant},
the throughput \eqref{othruput1} is the standard result  \eqref{globalthru},
and \eqref{othruput} follows from \eqref{localthru} and \eqref{trafficsolution}.
\end{beweis}
\begin{beweis} (\underline{of Corollary \ref{normingconstant-det}})
The simplest way is to approximate the general travel time distribution by a finite mixture of
Erlangian distributions. These mixtures constitute a class which is dense in the set
of all service times on
$[0,\infty)$, see \cite{daduna:01a}[Definition 9.2] and the references given there.
A Markovian state description of the network is obtained using supplementary variables.
Writing down the steady state throughput, we see that after some computations this
boils down to the same explicit expression as it occurs when writing down the throughput
expression for the companion exponential  network.\\
The last step is a continuity argument: When a sequence of finite mixtures of Erlangian service time distributions approaches (in the sense of weak convergence)
the given service time distribution at some node, then  that node's queue length distributions converge weakly as well.
This is equivalent to (multi-dimensional) point-wise convergence in a suitable multi-dimensional
real space. Finally, the normalization constants and the throughput continuously depend on the
 densities of the joint queue length distributions, as visible from Lemma
\ref{lem:normingconstant}.
\end{beweis}

\begin{lem}\label{lem:ungleichung}
For all $k=0,1,\ldots,N-2$, $n=0,1,\ldots,N-k-2$, it holds that
    \begin{displaymath}
        C_{n+k+1}(N,\bar{J})C_k(N-1,\bar{J})\ge C_{n+k+1}(N-1,\bar{J})C_k(N,\bar{J}).
    \end{displaymath}
\end{lem}
\begin{beweis}
    By definition, we have
    \begin{displaymath}
        C_k(N,\bar{J})=\left(\sum_{n_1+\ldots+n_J=N-k}\prod_{j=1}^J
        \left(\prod_{k=1}^{n_j}\frac{\eta_j}{\mu_j(k)}\right)^{n_j}\right)\frac{2^k}{k!},
    \end{displaymath}
   and hence with Theorem \ref{thm:gordonnewell}, $C_k(N,\bar{J})$ is the normalization constant $G(N-k,J)$ of a standard Gordon-Newell
    network with $N-k$ customers, $J$ service stations and not normalized solution $(\eta_1,\ldots,\eta_J)$
    of the associated traffic equation, multiplied with $2^k/k!$.
With these expressions and $\sum_{j=1}^J\eta_j=1/2$, it is easily verified that
    the throughput of this  Gordon-Newell network with $N$
    customers is $TH(N)=(1/2)G(N-1,J)/G(N,J)$.
    We then obtain
    \begin{align*}
        &\ C_{n+k+1}(N,\bar{J})C_k(N-1,\bar{J})\ge C_{n+k+1}(N-1,\bar{J})C_k(N,\bar{J})\\
        \Leftrightarrow &\ G(N-n-k-1,J)G(N-k-1,J)\ge G(N-n-k-2,J)G(N-k,J)\\
        \Leftrightarrow &\ \frac{G(N-k-1,J)}{G(N-k,J)}\ge\frac{G(N-n-k-2,J)}{G(N-n-k-1,J)}\\
        \Leftrightarrow &\  TH(N-k)\ge TH(N-n-k-1).
    \end{align*}
From \cite{vanderwal:89}, the throughput of a Gordon-Newell-network with service rates non-decreasing in the number of customers is a non-decreasing function in the network's population size. Thus, the lemma is proved.
\end{beweis}
\begin{beweis} \underline{(of Theorem \ref{thm:loesung})}
    For all $x\in\R^2$ with $h(x)>h(x^\ast)$ and $N\in\N$ we will show
    \begin{displaymath}
        TH(N;x^\ast)=\frac{G(N-1,\bar{J};x^\ast)}{G(N,\bar{J};x^\ast)}
        >\frac{G(N-1,\bar{J};x)}{G(N,\bar{J};x)}=TH(N;x)
    \end{displaymath}
   By Lemma \ref{lem:normingconstant} ,this is equivalent to
    \begin{eqnarray*}
        &   & G(N-1,\bar{J};x^\ast)G(N,\bar{J};x)-G(N-1,\bar{J};x) G(N,\bar{J};x^\ast)>0\\
        \Leftrightarrow &  & \left(\sum_{n=0}^{N-1} C_n(\N-1,\bar{J})h(x^\ast)^n)\right)
        \left(\sum_{k=0}^N       C_k(N,\bar{J})h(x)^k\right)\\
        & - & \left(\sum_{n=0}^{N-1} C_{n}(N-1,\bar{J})h(x)^n\right)
            \left(\sum_{k=0}^N C_k(N,\bar{J})h(x^\ast)^k)\right)>0\\
        \Leftrightarrow &  &
            \sum_{k=0}^N\sum_{n=0}^{N-1} C_k(N,\bar{J})
            C_n(N-1,\bar{J})\left(h(x)^kh(x^\ast)^n-h(x)^nh(x^\ast)^k\right)>0.
    \end{eqnarray*}
    We consider the summand for $k=N$
    \begin{align*}
        &    \sum_{n=0}^{N-1}C_N(N,\bar{J})C_n(N-1,\bar{J})
            \left(h(x)^Nh(x^\ast)^n-h(x)^nh(x^\ast)^N\right)\\
        & =  \sum_{n=0}^{N-1}C_N(N,\bar{J})C_n(N-1,\bar{J})
            h(x)^nh(x^\ast)^n\left(h(x)^{N-n}-h(x^\ast)^{N-n}\right).
    \end{align*}
    Because of $h(x)>h(x^\ast)$ we have $h(x)^{N-n}-h(x^\ast)^{N-n}>0$. So the whole summand is strictly positive
    and the problem is reduced to prove
    \begin{displaymath}
        \sum_{k=0}^{N-1}\sum_{n=0}^{N-1}C_{k}(N,\bar{J})
            C_{n}(N-1,\bar{J})(h(x)^kh(x^\ast)^n-h(x)^nh(x^\ast)^k)>0.
    \end{displaymath}
    For $k=n$ we get $h(x)^kh(x^\ast)^n-h(x)^nh(x^\ast)^k=0$. So the problem is  reduced to
    \begin{align*}
        &   \sum_{k=0}^{N-1}\sum_{n=0 \atop n\ne k}^{N-1}
            C_k(N,\bar{J})C_n(N-1,\bar{J})(h(x)^k
            h(x^\ast)^n-h(x)^nh(x^\ast)^k)>0\\
        \Leftrightarrow &
            \sum_{k=0}^{N-1}\sum_{n=k+1}^{N-1}
            C_k(N,\bar{J})C_n(N-1,\bar{J})(h(x)^k
            h(x^\ast)^n-h(x)^nh(x^\ast)^k)+\\
        &    \sum_{k=0}^{N-1}\sum_{n=0}^{k-1}
            C_k(N,\bar{J})C_n(N-1,\bar{J})(h(x)^k
            h(x^\ast)^n-h(x)^nh(x^\ast)^k)>0.
    \end{align*}
    Now consider the case $k=N-1$ in the first summand. The second sum is empty ($=0$). The same holds in the second
    summand for $k=0$. So we have reduced the problem  to
    \begin{eqnarray*}
        &   & \sum_{k=0}^{N-2}\sum_{n=k+1}^{N-1}
            C_k(N,\bar{J})C_n(N-1,\bar{J})(h(x)^k
            h(x^\ast)^n-h(x)^nh(x^\ast)^k)\\
        & + & \sum_{k=1}^{N-1}\sum_{n=0}^{k-1}
            C_k(N,\bar{J})C_n(N-1,\bar{J})(h(x)^k
            h(x^\ast)^n-h(x)^nh(x^\ast)^k)>0.
    \end{eqnarray*}
    by index-shift in both summands we get
    \begin{eqnarray*}
        &   & \sum_{k=0}^{N-2}\sum_{n=0}^{N-k-2}
            C_k(N,\bar{J})C_{n+k+1}(N-1,\bar{J})(h(x)^k
            h(x^\ast)^{n+k+1}-h(x)^{n+k+1}h(x^\ast)^k)\\
        & + & \sum_{k=0}^{N-2}\sum_{n=0}^{k}
            C_{k+1}(N,\bar{J})C_n(N-1,\bar{J})(h(x)^{k+1}
            h(x^\ast)^n-h(x)^nh(x^\ast)^{k+1})>0.
    \end{eqnarray*}
    In the second summand we apply  the following summation formula twice:
    \begin{displaymath}
         \text{for}~~~a_k,b_k\in\R,~~k=1,\ldots,N,~~~~
        \sum_{k=0}^N\sum_{n=0}^k a_{k+1}b_n
        =\sum_{k=0}^N\sum_{n=0}^{N-k} a_{n+k+1}b_k,
    \end{displaymath}
   and obtain
    \begin{eqnarray*}
        &   & \sum_{k=0}^{N-2}\sum_{n=0}^{N-k-2}
            C_k(N,\bar{J})C_{n+k+1}(N-1,\bar{J})(h(x)^k
            h(x^\ast)^{n+k+1}-h(x)^{n+k+1}h(x^\ast)^k)\\
        & + & \sum_{k=0}^{N-2}\sum_{n=0}^{N-k-2}
            C_{n+k+1}(N,\bar{J})C_k(N-1,\bar{J})(h(x)^{n+k+1}
            h(x^\ast)^k-h(x)^kh(x^\ast)^{n+k+1})>0.
    \end{eqnarray*}
    Because of $h(x)>h(x^\ast)$, we have $h(x)^{n+k+1}h(x^\ast)^k-h(x)^kh(x^\ast)^{n+k+1}>0$
    and from Lemma \ref{lem:ungleichung} we have $C_{n+k+1}(N,\bar{J})C_k
    (N-1,\bar{J})\ge C_{n+k+1}(N-1,\bar{J})C_k(N,\bar{J})$. So
    \begin{eqnarray*}
        &   & \sum_{k=0}^{N-2}\sum_{n=0}^{N-k-2}
            C_k(N,\bar{J})C_{n+k+1}(N-1,\bar{J})(h(x)^k
            h(x^\ast)^{n+k+1}-h(x)^{n+k+1}h(x^\ast)^k)\\
        & + & \sum_{k=0}^{N-2}\sum_{n=0}^{N-k-2}
            C_{n+k+1}(N,\bar{J})C_k(N-1,\bar{J})(h(x)^{n+k+1}
            h(x^\ast)^k-h(x)^kh(x^\ast)^{n+k+1})\\
        &\ge &\sum_{k=0}^{N-2}\sum_{n=0}^{N-k-2}
            C_k(N,\bar{J})C_{n+k+1}(N-1,\bar{J})(h(x)^k
            h(x^\ast)^{n+k+1}-h(x)^{n+k+1}h(x^\ast)^k)\\
        & + & \sum_{k=0}^{N-2}\sum_{n=0}^{N-k-2}
            C_k(N,\bar{J})C_{n+k+1}(N-1,\bar{J})(h(x)^{n+k+1}
            h(x^\ast)^k-h(x)^k h(x^\ast)^{n+k+1})=0,
    \end{eqnarray*}
    and  the theorem is  proved.
\end{beweis}

\subsection{Additional numerical experiments}\label{sect:AddNumeric}
We report in this section  details of the  experiments  which have been summarized in Table \ref{tab:Summary} in  Section \ref{sect:Numeric}. We performed four blocks
of experiments, distinguished by  demand distributions which
are uniform $U(A)$ on  finite demand sets $A$.\\
In any case  $12$ locations are sampled uniformly from
$\{10,\dots,410\}\times \{10,\dots,270\}$.
Service is provided by single servers with intensities
 $\mu_j=2, j=2,\dots,13$, for unloading servers at warehouses.
Service intensity $\mu_1$ at the center varies with the blocks.\\  
\underline{Abbreviations:}\\  
 DistLoc $\equiv$ distance between center  with weights and center  without weights\\
Demand: to/$\wedge$ /$\vee$ $\equiv$ total demand/minimal demand/maximal demand \\
+/- $\equiv$  quantities for: center with weights(+)/ center without weights (-)

The extreme cases of large demand are highlighted by boldface numbers, the extreme differences for needed trucks are colored red.
The cases where the requested demand could not be delivered due to bounds determined by bottlenecks are indicated under ``Trucks''
as ``--/--''. Bottleneck was in any case the loading server at the center.
{\small
\begin{table}
		\begin{tabular}{|r|r|r|r|r|}\hline
		DistLoc&Demand&Trucks&Throughput/day& $P(X_1>0)$ \\ \hline
		& to/$\wedge$/$\vee$  & +/-~~  &          +/-~~~~~~~~   &    +/- ~~~~   \\ \hline
		19.0874& 52/1/8& 12/12& 53.4447/53.1933& .5567/.5541  \\ \hline 
		11.2037& 54/2/8& 13/13& 56.5579/56.2689& .5891/.5861 
		\\ \hline   
		14.4548& 50/1/8& 11/11& 52.1731/51.9762& .5435/.5414
		\\ \hline   
		\textbf{60.7717}& 41/1/8& 9/10& 42.5902/43.5175& .4436/.4533
		\\ \hline   
		7.5896& 44/1/7& 10/10& 46.1790/46.1360& .4810/.4806
		\\ \hline   
		36.9812& 59/1/8& 16/17& 59.5479/61.9236& .6203/.6450
		\\ \hline   
		52.2486& 49/1/8& 11/12& 50.6667/52.9925& .5278/.5520
		\\ \hline   
		23.0776& 59/2/8& 14/14& 61.5728/61.0184& .6414/.6356
		\\ \hline   
		53.3149& 45/1/7& 13/14& 45.8186/47.9310& .4773/.4993
		\\ \hline  
		\textbf{57.3881}& 63/1/8& 19/19& 64.6827/64.0564& .6738/.6673
		\\ \hline   
	\end{tabular} 
\caption{\textbf{(I)}:  Demand distribution $U(\{1,\dots,8\})$, loading rate  $\mu_1=4$.}
\end{table}

\begin{table}
	\begin{tabular}{|r|r|r|r|r|}\hline
		DistLoc&Demand&Trucks&Throughput/day& $P(X_1>0)$ \\ \hline
		& to/$\wedge$/$\vee$  & +/-~~  &          +/- ~~~~~~~~~  &    +/-~~~ ~~   \\ \hline
		\textbf{88.8410}& 107/2/16& 29/31& 107.1006/108.5612& .8925/.9047
		\\ \hline  
		29.8252& 111/2/16& 32/32& 112.6342/111.6704& .9386/.9306
		\\ \hline   
		42.3096& 112/2/16& 34/34& 113.1859/112.4563& .9432/.9371
		\\ \hline  
		54.8994& 98/2/15& 25/26& 98.0313/98.4849& .8169/.8207
		\\ \hline   
		16.5874& 121/3/16& --/--& 120.0000/120.0000& 1.0000/1.0000
		\\ \hline  
		24.7907& 93/2/15& 24/24& 95.3214/94.5225& .7943/.7877
		\\ \hline   
		\textbf{76.2196}	& 96/2/16& 27/28& 96.6996/96.7557& .8058/.8063
		\\ \hline  
		4.9971& 123/2/16& --/--& 120.0000/119.9999& 1.0000/.9999
		\\ \hline   
		41.8285& 97/3/16& 27/28& 97.8266/98.2422& .8152/.8187
		\\ \hline  
		16.2628& 105/3/16& 29/29& 106.3469/105.8383& .8862/.8820
		\\ \hline   
		51.2422& 99/1/16& 28/29& 99.2529/100.0747& .8271/.8340
		\\ \hline  
		30.2271& 107/1/16& 31/32& 107.1240/107.9508& .8927/.8996
						\\ \hline  
	\end{tabular}
\caption{\textbf{(II)}:  Demand distribution $U(\{1,\dots,16\})$, loading rate $\mu_1=5$.}
\end{table}

\begin{table}
\begin{tabular}{|r|r|r|r|r|}\hline
	DistLoc&Dem&Tru&Throughput/day& $P(X_1>0)$ \\ \hline
	& to/$\wedge$/$\vee$  & +/-~~  &          +/-~~~~~~~~~~   &    +/-~~~~~    \\ \hline
	16.8426& 163/3/20& 43/44& 163.2859/163.9906& .9719/.9761
	\\ \hline   
	28.9316& 138/2/19& 37/38& 138.9829/140.1825& .8273/.8344
	\\ \hline  
	\textbf{108.5533}& 101/2/21& 
	\underline{27/30}&
	102.7856/101.7569& .6118/.6057
	\\ \hline   
	19.0861& 113/1/18& 28/28& 113.9980/113.5362& .6786/.6758
	\\ \hline  
	6.6153& 188/6/21& --/--& 167.9999/168.0000& .9999/1.0000	
	\\ \hline   
	11.0563& 127/4/18& 29/29& 129.1835/128.7957& .7689/.7666
	\\ \hline  	
	24.5071& 102/1/21& 24/24& 105.0053/103.6437& .6250/.6169
	\\ \hline   
	26.8209& 94/1/20& 27/27& 96.2469/94.8691& .5729/.5647
	\\ \hline  
	18.4939& 156/1/21&45/46& 156.3308/157.6311& .9305/.9383
	\\ \hline   
	16.6062& 156/3/21& 32/32& 157.5397/156.5776& .9377/.9320
	\\ \hline  
	\textbf{39.9891}& 119/2/21& 31/32& 120.2597/120.0498& .7158/.7146\\ \hline
\end{tabular}
\caption{\textbf{(III)}:  Demand distribution $U(\{1,\dots,21\})$, loading rate $\mu_1=7$.}
\end{table}

\begin{table}
\begin{tabular}{|r|r|r|r|r|}\hline
	DistLoc&Dem&Tru&Throughput/day& $P(X_1>0)$ \\ \hline
	& to/$\wedge$/$\vee$  & +/-~~  &          +/-~~~~~~~~~~   &    +/-~~~~~    \\ \hline
	35.0573& 142/1/21& 37/37& 144.2880/142.0597& .8589/.8456
	\\ \hline  
	\textbf{85.6094}& 122/1/21& 
	\underline{33/36}&
	124.7843/123.4753& .7428/.7350
	\\ \hline   
	29.3856& 132/1/21& 32/33& 132.0441/134.0093& .7860/.7977
	\\ \hline  	
	.1719& 202/11/21& --/--& 167.9999/168.0000& .9999/1.0000
	\\ \hline   
	21.9623& 162/11/21& 54/55& 162.0168/162.8296& .9644/.9692
	\\ \hline  
	5.5901& 142/1/21& 31/31& 144.8758/144.8089& .8624/.8620
	\\ \hline   
	37.6852& 42/1/11& 6/7& 42.6926/44.3910& .2541/.2642
	\\ \hline  
	18.3488& 72/1/21& 17/17& 75.1295/74.8134& .4472/.4453
	\\ \hline   
	5.5115& 152/1/21& 38/38& 152.1765/152.1388& .9058/.9056
	\\ \hline  
	6.0370& 192/11/21&--/--& 168.0000/ 167.9999& 1.0000/.9999
	\\ \hline   
	\textbf{48.1097}& 132/1/21& 33/34& 134.0373/133.2210& .7978/.7930	
	\\ \hline  	
	21.2497& 122/1/21& 30/31& 122.7658/124.9653& .7307/.7438
	\\ \hline   
\end{tabular}
\caption{\textbf{(IV)}:  Demand distribution $U(\{1,11,21\})$, loading rate $\mu_1=7$.}
\end{table}
}
\end{appendix}

\textbf{Acknowledgment:} We thank Peter Sieb for helpful discussions on the subject of this paper.

\end{document}